\documentclass[11pt]{article}
\usepackage[english]{babel}
\usepackage{geometry}
\usepackage{setspace}
\usepackage{graphicx}
\graphicspath{{.}}
\usepackage{lineno}
\usepackage{multirow}
\usepackage{babel}
\usepackage{booktabs}
\usepackage{footnote}
\usepackage{threeparttable}
\usepackage{color}

\usepackage{floatrow}
\newfloatcommand{capbtabbox}{table}[][\FBwidth]

\usepackage{amsmath,amsthm, amssymb}

\usepackage{color}
\usepackage[toc,page]{appendix}
\usepackage[colorlinks,linkcolor=black,anchorcolor=blue,citecolor=blue,CJKbookmarks=True]{hyperref}

\usepackage[T1]{fontenc}
\usepackage[utf8]{inputenc}
\usepackage{babel}
\usepackage[font=small,labelfont=bf,tableposition=top]{caption}
\usepackage{booktabs}
\usepackage{threeparttable}
\usepackage{footnote}
\usepackage{multirow}
\usepackage{natbib}
\bibpunct[, ]{(}{)}{,}{a}{}{,}%

\usepackage{algorithm}
\usepackage{algpseudocode}
\usepackage{eqparbox}

\title{
Joint optimization of train blocking and shipment path:\\
An integrated model and a sequential algorithm}

\author{Chongshuang Chen$^{1,2}$,\quad Jun Zhao$^3$ \\
ccsmars@swjtu.edu.cn,\quad junzhao@swjtu.edu.cn\\[0.5cm]
\small{$^1$ School of Mathematics, Southwest Jiaotong University, Chengdu, Sichuan 611756, China}\\
\small{$^2$ National Engineering Laboratory of Integrated Transportation Big Data Application Technology,}\\
\small{Chengdu, Sichuan 611756, China}\\
\small{$^3$ School of Transportation and Logistics, Southwest Jiaotong University, Chengdu, Sichuan 611756, China}}
\date{}

\begin{document}
\maketitle

\section{Introduction}
\label{sec: introduction}
The INFORMS RAS 2019 Problem Solving Competition is focused on the integrated train blocking and shipment path (TBSP) optimization for tonnage-based operating railways. In nature, the TBSP problem could be viewed as a multi-commodity network design problem with a double-layer network structure. By introducing a directed physical railway network and a directed train services (blocks) network, we formulate completely the TBSP problem as a mixed integer linear programming (MILP) model that incorporates all decisions, objectives and constraints especially the merge flow (called intree rule here) in an integrated manner. The scale of the MILP model can be reduced efficiently if we only enumerate the arc selection and block sequence variables for each shipment on the legal paths from its origin to its destination satisfying the given detour ratio. We further develop a sequential algorithm that decomposes the TBSP problem into the shipment path subproblem and train blocking subproblem which are solved sequentially. Computational tests on the three given data sets show that the reduced MILP model can solve DataSet\_1 to optimality in 8.48 seconds and DataSet\_2 with a gap of 0.16\% in 6 hours on a GPU workstation. The reduced model also can provide strong lower bounds for DataSet\_2 and DataSet\_3. The sequential algorithm can find a high quality solution with 0.04\% gap within 0.26 seconds for DataSet\_1, 0.42\% gap within 4.53 seconds for DataSet\_2 and 1.54\% gap within 0.58 hours for DataSet\_3 respectively on a Thinkpad laptop.
%The remainder of this paper is organized as follows. An integrated formulation for the TBSP problem and model reduction technologies are introduced in \S\ref{sec: Integrated formulation}. The sequential formulation is presented in \S\ref{sec: Sequential formulation}. \S\ref{sec: Computational Results} discusses computational results, and \S\ref{sec: Conclusions} summarizes our conclusions.

\section{Integrated model}
\label{sec: Integrated formulation}
In this section, we develop an integrated MILP model for the TBSP problem by introducing two directed graphs. The aim is to make decisions on the shipment path, train services (blocks) design, train frequency and shipment-block sequence (car-to-block assignment). The objective is the minimization of the total sum of the car transportation cost, train accumulation delay and car reclassification delay, subject to basic network design constraints and the intree rule.

\subsection{Modelling motivation}
\label{sec: Approach motivation}
By definition, the TBSP problem has two interconnected subproblems including the shipment path and the train blocking. The former identifies the path of shipments to minimize the total car transportation cost and respect the capacity of railway lines. The latter determines the train services, train frequencies and shipment-block sequences such that the classification capacity and sort tracks in yards as well as the intree rule are satisfied, while the total train accumulation delay and car reclassification delay are minimized. Under an integrated framework, there are intricate relationships between these two subproblems \citep{Lin2012}. Specifically, every shipment must be transported from its origin to its destination on the physical railway network. Meanwhile, all shipments must be moved by the means of train services (blocks) with necessary accumulation and/or reclassification operations in the origin yard, intermediate yard(s) and destination yard. Importantly, the reclassification yard(s) of each shipment should be limited on its path, and the path of a train service must be completely consistent with that of the shipments transported on the train service.

\begin{figure}[h!]
  \caption{A simple network}\label{fig: A simple network}
  \centering
  \includegraphics[scale=0.8]{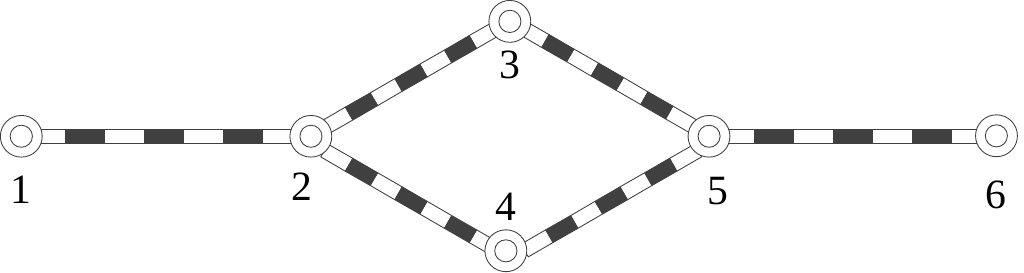}
\end{figure}

\begin{figure}[h!]
  \caption{Double-layer network structure}\label{fig: Double layer network structure}
  \centering
  \includegraphics[scale=0.62]{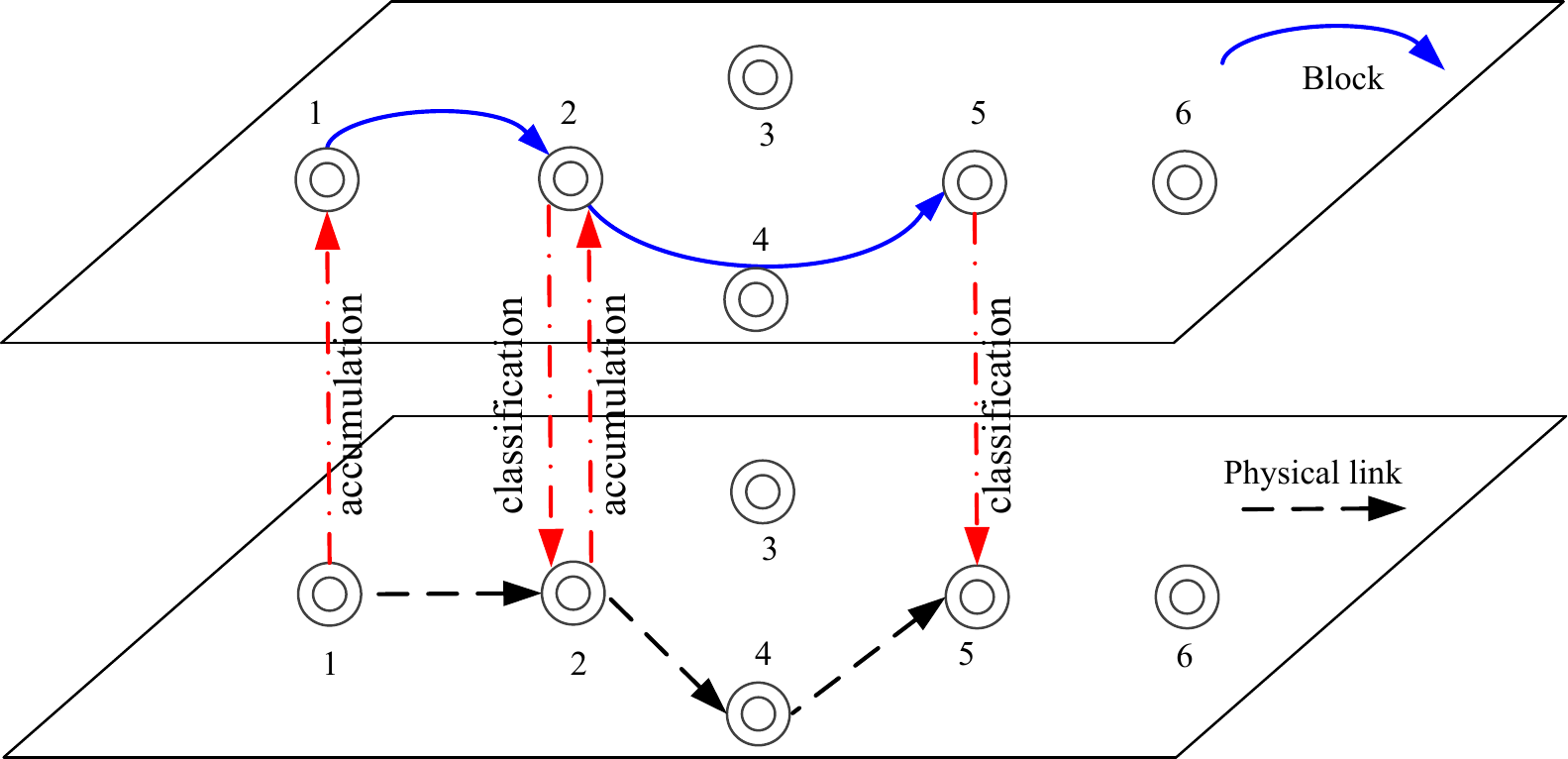}
\end{figure}

We present an example with 6 yards and 6 links in Figure \ref{fig: A simple network}
to explain. For the shipment from 1 to 5, we list all its possible paths, block sequences and reclassification yards in Table \ref{tab: Relationship between paths and block sequences}. Assume that this shipment chooses path $1\rightarrow2\rightarrow4\rightarrow5$ and block sequences $1\rightarrow2$, $2\rightarrow5$, respectively. Then, the reclassification yards of the shipment can only be 2 and 5. Meanwhile, the path decision of the train service from $2$ to $5$ should be $2\rightarrow4\rightarrow5$. Figure \ref{fig: Double layer network structure} roughly shows the transportation process for the shipment from $1$ to $5$ on the network in Figure \ref{fig: A simple network}. In nature, the TBSP problem could be viewed as a multi-commodity network design problem with a double-layer network structure.

\begin{table}[h!]\small
\begin{center}
\caption{Relationship between shipment path, block sequences and reclassification yards}
\label{tab: Relationship between paths and block sequences}
\begin{tabular}{ c l l}
\hline
Shipment path &  Block sequences & Reclassification yards\\
\hline
\multirow{4}{*}{$1\rightarrow2\rightarrow3\rightarrow5$}
& $1\rightarrow5$      &5\\
& $1\rightarrow2$, $2\rightarrow5$  & 2, 5\\
& $1\rightarrow3$, $3\rightarrow5$  & 3, 5\\
& $1\rightarrow2$, $2\rightarrow3$, $3\rightarrow5$  & 2, 3, 5\\
\hline
\multirow{4}{*}{$1\rightarrow2\rightarrow4\rightarrow5$}
& $1\rightarrow5$      &5\\
& $1\rightarrow2$, $2\rightarrow5$  & 2, 5\\
& $1\rightarrow4$, $4\rightarrow5$  & 4, 5\\
& $1\rightarrow2$, $2\rightarrow4$, $4\rightarrow5$  & 2, 4, 5\\
\hline
\end{tabular}
\end{center}
\end{table}

We introduce two directed networks and define two arc-based decision variables to route all shipments through the double-layer network.
We also utilize well-designed constraints to reflect the above discussed interconnections between the shipment path and train blocking subproblems of the TBSP problem in our integrated model and sequential algorithm.
For more details, please see \S\ref{sec: MIQP formulation} and \S\ref{sec: Blocking problem}.

\subsection{Parameters and variables}
\label{sec: Parameters and variables}
The problem is referred to two directed graphs. One is the physical railway network $G_1=(V,E)$, where $V$ denotes the set of all yards, and $E \subseteq V\times V$ the set of links between two adjacent yards. Another is the service network $G_2=(V,\mathcal{B})$, where $\mathcal{B} = V\times V$ is the set of all potential blocks (directed train services), i.e., one for any two yards. The latter is a complete graph which is much denser than the former.

We list the parameters of the TBSP problem in Table \ref{tab: parameter}. Table \ref{tab: Decision Variables} lists all kinds of decision variables used in our integrated model. To be greatly different from classical service network design problems, a compulsive regulation arises in the TBSP problem is called the \emph{intree rule} in this report. Such a rule states that all shipments with the same destination must be operated in an exact same way once they are reclassified at the same intermediate yard.
%Thus, when a demand arrives at an intermediate yard, the yard it must be sent to next is regardless of the demand's origin.
We follow the scheme applied by \cite{Powell1992} to deal with the tree constraint in less-than-truckload transportation. Except four decisions required in the problem description, another binary variable ($v^{d}_{p,q}$) is introduced to ensure the intree rule to be hold. The last one addresses the sort track constraints.
\begin{table}[h!]\small
\begin{center}
\caption{Parameters of the TBSP}
\label{tab: parameter}
\begin{tabular}{ c l l l}
\hline  Notation & Definition & Unit & Remark  \\
\hline $V$ & set of yards & & index $o,d,p,q,i,j $ \\
$E$ & set of links &  &index $e=(i,j) $ \\
%$A_{i}^{+}$ & set of successors of yard $i$, $\forall i\in V$ &   & $A_{i}^{+}=\{j \in V |(i,j)\in E \}$\\
%$A_{i}^{-}$ & set of predecessors of yard $i$, $\forall i\in V$ &   & $A_{i}^{-}=\{j \in V |(j,i)\in E \}$\\
$\mathcal{B}$ & set of blocks & & index $B_{p,q}$\\
$l_{i,j}$ & length of link  $(i,j)$, $ \forall (i,j)\in E$ &  km & constant \\
$t_i$ & reclassification delay time per car in yard $i$, $\forall i \in V$ & hour & constant\\
$c_{p,q}$ & accumulation parameter for train from $p$ to $q$, $\forall p,q \in V$ & hour & constant\\
$m$ & train size  &  car & constant\\
$n_{o,d}$ & cars of shipment from $o$ to $d$, $\forall o,d \in V$  &  car & constant \\
$f_{i,j}$ & capacity on link $(i,j)$, $ \forall (i,j)\in E$ &  train & constant \\
$g_i$ & original classification capacity in yard $i$, $ \forall i\in V$ &  car & constant \\
$h_i$ & number of sort tracks in yard $i$, $ \forall i\in V$ & & constant \\
$\alpha_{i,j}$ & remaining rate of capacity on link $(i,j)$, $ \forall (i,j)\in E$ &  & constant \\
$\beta_i$ & available ratio of capacity in yard $i$, $ \forall i\in V$ & & constant \\
$\gamma$ & reclassification capacity of each sort track & car & constant \\
$\delta_{o,d}$ & length of the shortest path from $o$ to $d$, $ \forall o,d\in V$ &  km & constant \\
$\epsilon$ & detour ratio threshold    &   & constant \\
$\lambda$ & conversion factor for car-kilometers   &   & constant \\
\hline
\end{tabular}
\end{center}
\end{table}

\begin{table}[h!]\small
\begin{center}
\caption{Decision variables of the TBSP}
\label{tab: Decision Variables}
\begin{tabular}{ c l l }
\hline  Notation & Definition & Remark\\
\hline
$x^{o,d}_{i,j}$ & 1 if shipment(block) from $o$ to $d$ passes through link $(i,j)$, 0 otherwise & arc selection \\
$y_{p,q}$ & 1 if block $B_{p,q}$ is provided, 0 otherwise & block design\\
$z_{p,q}$ & frequency of train carries block $B_{p,q}$  & train frequency\\
$u^{o,d}_{p,q}$ & 1 if shipment from $o$ to $d$ is consolidated into block $B_{p,q}$, 0 otherwise & block sequence \\
$v^{d}_{p,q}$ & 1 if shipment from $p$ to $d$ is consolidated into block $B_{p,q}$, 0 otherwise & consolidation selection \\
$w_{p,q}$ & number of sort tacks needed by block $B_{p,q}$  & sort tack usage\\
\hline
\end{tabular}
\end{center}
\end{table}

\subsection{MIQP model}
\label{sec: MIQP formulation}
Using above parameters and decision variables, the integrated train blocking and shipment path problem is modelled as a mixed integer quadratic programming (MIQP) as follows.
\begin{equation}\small
\label{eq: objective function}
\text{(MIQP)} \quad \min \quad \lambda\sum_{o,d\in V}\sum_{(i,j)\in E}n_{o,d}l_{i,j}x_{i,j}^{o,d}+\sum_{p,q\in V}mc_{p,q}y_{p,q}
+ \sum_{o,d\in V}\sum_{\{p,q \in V | q\neq d\}}n_{o,d}t_qu^{o,d}_{p,q}
\end{equation}
\begin{align}\small
&\text{s.t.}&\sum_{ \{j \in V |(i,j)\in E \}}x_{i,j}^{o,d}-\sum_{ \{j \in V |(j,i)\in E \}}x_{j,i}^{o,d}&
=\left\{\begin{array}{ll}
             1, &\mbox{if $i=o$} \\
            -1, & \mbox{if $i=d$} \\
             0, & \mbox{if $i \neq o,d$} \\
\end{array}\right. &&\forall o,d\in V \label{constraint: car balance}\\
&&\sum_{p,q\in V }x_{i,j}^{p,q}z_{p,q} &\leq f_{i,j}\alpha_{i,j} &&\forall (i,j) \in E \label{constraint: link capacity}\\
&&\sum_{(i,j)\in E}l_{i,j}x_{i,j}^{o,d} &\leq \epsilon\delta_{o,d} &&\forall o,d \in V \label{constraint: detour}\\
&&\sum_{o,d\in V }n_{o,d}u^{o,d}_{p,q}-mz_{p,q} & =0&&\forall p,q \in V \label{constraint: transport capacity} \\
&&\sum_{q\in V }u^{o,d}_{p,q}-\sum_{q\in V}u^{o,d}_{q,p}&
=\left\{\begin{array}{ll}
             1, &\mbox{if $p=o$} \\
            -1, & \mbox{if $p=d$} \\
             0, & \mbox{if $p \neq o,d$} \\
\end{array}\right. &&\forall o,d\in V \label{constraint: block balance}\\
&&\sum_{p \in V} \sum_{\{o,d \in V | d\neq q\}}n_{o,d}u^{o,d}_{p,q} &\leq g_q\beta_q &&\forall q \in V \label{constraint: yard capacity}\\
&&\sum_{q \in V} w_{p,q} &\leq h_p && \forall p\in V \label{constraint: yard track}\\
&& \gamma (w_{p,q}-1)+1 &\leq \sum_{o,d\in V }n_{o,d}u^{o,d}_{p,q} &&\forall p,q \in V \label{constraint: track usage lower bound}\\
&& \sum_{o,d\in V }n_{o,d}u^{o,d}_{p,q}-\gamma w_{p,q} &\leq 0&&\forall p,q \in V \label{constraint: track usage upper bound}\\
&&u_{p,q}^{o,d}-v_{p,q}^d &\leq 0 &&\forall o, d,p,q \in V \label{constraint: exclusive rule 1}\\
&&\sum_{q\in V }v_{p,q}^d &\leq 1 &&\forall d,p \in V \label{constraint: exclusive rule 2} \\
&&u_{p,q}^{o,d}-y_{p,q} &\leq 0 &&\forall o, d,p,q \in V \label{constraint: car-blcok cross}\\
&&u_{p,q}^{o,d}+x_{i,j}^{p,q}-x_{i,j}^{o,d} &\leq 1 &&\forall o, d,p,q \in V,(i,j) \in E \label{constraint: path consistant}
\end{align}
\begin{align}\small
&&x^{o,d}_{i,j} &\in \{0,1\} &&\forall o,d \in V, (i,j) \in E \label{constraint: arc selection variables domain}\\
&&y_{p,q} &\in \{0,1\} &&\forall p,q \in V \label{constraint: block design variables domain}\\
&&z_{p,q} &\geq 0 &&\forall p,q \in V \label{constraint: train frequency variables domain}\\
&&u_{p,q}^{o,d} &\in \{0,1\} &&\forall o,d,p,q \in V \label{constraint: car distribution variables domain} \\
&&v_{p,q}^d &\in \{0,1\} &&\forall d,p,q \in V \label{constraint: exclusive rule variables domain}\\
&&w_{p,q} &\in \mathbb{Z} &&\forall p,q \in V \label{constraint: track usage variables domain}
\end{align}

The objective function \eqref{eq: objective function} is to minimize the total sum of car transportation cost, train accumulation delay and car reclassification delay at intermediate stops in shipping all demands over an underlying network. The first part is measured by car-kilometer whereas the last two are measured by car-hour. Thus, a conversion factor is introduced to balance off these two different units.

Flow conservation equations \eqref{constraint: car balance} guarantee that each shipment can reach its destination on the physical network. Inequalities \eqref{constraint: link capacity} ensure the number of trains passing on a link does not exceed its capacity. Constraints \eqref{constraint: detour} are presented to prevent each shipment from transporting through a too long detour path. Tight constraints \eqref{constraint: transport capacity} calculate the train frequency of each provided block.

Flow conservation equations \eqref{constraint: block balance} guarantee that every shipment can reach its destination on the service network.
%All blocks with the same destination will share their common arrival yard's capacity.
Inequalities \eqref{constraint: yard capacity} ensure that the amount of cars reclassified at each yard does not exceed the available capacity of the yard.
%It should be noted that all demands would be definitely classified at their destinations no matter a direct-train or some connected-trains are adopted. Therefore, the objective function \eqref{eq: objective function} does not include the delay cost at the destination, while the constraint \eqref{constraint: yard capacity} dose not consider the classification workload at the destination.
Inequalities \eqref{constraint: yard track} assure that the usage of sort tracks at each yard is under the budget of the available quantity. Constraints \eqref{constraint: track usage lower bound} and \eqref{constraint: track usage upper bound} indicate that total number of cars departing from each yard is strictly bounded by the sort track usage variables.
%Taking the summation with respect to $q$ in inequality \eqref{constraint: track usage upper bound} and combining with inequality \eqref{constraint: yard track}, we get an aggregated constraint $\sum_{q \in V}\sum_{o,d\in V }n_{o,d}u^{o,d}_{p,q} \leq \gamma h_p$. That is to say, aggregated constraints about sort track consumption are looser than disaggregated ones.

Disaggregate inequalities \eqref{constraint: exclusive rule 1} indicate the interior logical relationship between the block sequence variables and consolidation selection variables. Constraints \eqref{constraint: exclusive rule 2} ensure the uniqueness of the consolidation strategy. If shipments destined to yard $d$ are currently reclassified at yard $p$, then they must be merged into almost one block outward from $p$ to $d$. Both \eqref{constraint: exclusive rule 1} and \eqref{constraint: exclusive rule 2} enforce the intree rule to be respected. Constraints \eqref{constraint: car-blcok cross} forbid any shipment to be consolidated into the blocks that are not provided. Constraints \eqref{constraint: path consistant} trace the relationship between shipment paths and block sequence scenarios. To be specific, once shipment $o \rightarrow d$ is consolidated into block $B_{p,q}$, i.e., $u_{p,q}^{o,d}=1$, and link $(i,j)$ is on the route of $B_{p,q}$, i.e., $x_{i,j}^{p,q}=1$, then link $(i,j)$ must be on the route of demand $o \rightarrow d$, i.e., $x_{i,j}^{o,d}=1$.
%In other words, the route of the train (block) must be the intersections of demands who use it.
Constraints \eqref{constraint: arc selection variables domain}-\eqref{constraint: track usage variables domain} specify the domain of decision variables.

%Flow conservation equalities \eqref{constraint: car balance} and \eqref{constraint: block balance} are provided on the physical network and block network, respectively. In nature, the TBSP problem could be viewed as a multi-commodity flow problem with a double-layer network structure. It also should be noted that train frequency variables in \eqref{constraint: train frequency variables domain} are settled as non-negative continuous not integer. Because the volume of railcars are taken the representative quantity (e.g., average value) during the scheduling period. Therefore, both the train blocking and shipment path belong to tactical level decisions not operational level \citep{Cordeau1998}. The MIQP model could be extended to the situation when taking integer train frequency and different train sizes. A little change is tight equalities \eqref{constraint: transport capacity} would be slacked as inequalities.

\subsection{MILP model}
\label{sec: MILP formulation}
The MIQP model has a linear objective function and quadratic constraints \eqref{constraint: link capacity}.
%What we want to do is to try to linearize our nonlinear model, aiming at applying methodologies and algorithms in the integer linear programming. Constraint \eqref{constraint: link capacity} in
The MIQP model could be easily reconstructed by introducing another continuous variables $s^{p,q}_{i,j}=x^{p,q}_{i,j}z_{p,q}$ and auxiliary constraints \eqref{constraint: auxiliary constraint 1}-\eqref{constraint: auxiliary constraint 3}.
\begin{align} \small
&&s_{i,j}^{p,q}-Mx^{p,q}_{i,j} &\leq 0 &&\forall p,q\in V,(i,j) \in E \label{constraint: auxiliary constraint 1}\\
&&s_{i,j}^{p,q}-z_{p,q} &\leq 0 &&\forall  p,q\in V,(i,j) \in E \label{constraint: auxiliary constraint 2}\\
&&z_{p,q}+Mx^{p,q}_{i,j}-s_{i,j}^{p,q} &\leq M&&\forall p,q\in V,(i,j) \in E \label{constraint: auxiliary constraint 3}\\
&&s_{i,j}^{p,q} & \geq0 &&\forall p,q\in V,(i,j) \in E \label{constraint: auxiliary variables domain 1}
\end{align}
where $M$ in constraints \eqref{constraint: auxiliary constraint 1} and \eqref{constraint: auxiliary constraint 3} are big positive constants which could be set as $\frac{1}{m}\sum_{o,d \in V}n_{o,d}$.
Then, the TBSP problem is reformulated as a mixed integer linear programming (MILP) as follows.
\begin{equation*}\small
\text{(MILP)}\quad \min \eqref{eq: objective function}
\end{equation*}
\begin{align}\small
&\text{s.t.}&\sum_{p,q\in V }s_{i,j}^{p,q} &\leq f_{i,j}\alpha_{i,j} &&\forall (i,j)\in E \label{constraint: linear link capacity}
\end{align}
\begin{equation*}\small
  \eqref{constraint: car balance},\eqref{constraint: detour}-\eqref{constraint: auxiliary variables domain 1}.
\end{equation*}

The size of the MILP model can be roughly estimated as $|V|^4+|V|^3+2|V|^2|E|+3|V|^2$ variables and $|V|^4|E|+2|V|^4+3|V|^2|E|+2|V|^3+5|V|^2+2|V|$ constraints. Most variables are from the block sequence decisions ($u^{o,d}_{p,q}, \forall o,d,p,q \in V$), whereas most constraints arise from \eqref{constraint: path consistant} to keep the legality of all shipment-block sequences.

\subsection{Model reduction}
\label{sec: Model reduction}
As stated in the previous subsection, the scale of the MILP model increases rapidly as the problem size grows. We propose some model reduction technologies in this subsection. In Table \ref{tab: Decision Variables}, we set the arc selection variables ($x^{o,d}_{i,j}, \forall o,d \in V, (i,j)\in E$) and block sequence variables ($u^{o,d}_{p,q}, \forall o,d,p,q \in V$) using completed enumeration. We remark that such a manner is simple but not wise if the path detour ratio of shipments is enforced. We demonstrate the explanations as follows.
\begin{itemize}
  \item For each shipment, only a limited number of paths without too much detour could be allowed to be candidates. Due to that, not all links in the physical network $G_1$ could be passed over by each shipment. We denote $E_{o,d}$ as the set of links in $G_1$ that could be used by the shipment from $o$ to $d$.
  \item %If we view a path of a shipment as an order set along the direction from the origin to the destination of the shipment, then the block sequences could be regarded as an order subset of the path as well, as shown in Table \ref{tab: Relationship between paths and block sequences}.
      For each shipment, because of the restriction of path selection, not all blocks in the service network $G_2$ could be adopted by each shipment. We denote $\mathcal{B}_{o,d}$ as the set of candidate blocks in $G_2$ that could be available to the shipment from $o$ to $d$.
\end{itemize}

% if the defined detour ratio $\epsilon$ is bigger than 1.
We can enumerate variables $x^{o,d}_{i,j}$ on $E_{o,d}$ and variables $u^{o,d}_{p,q}$ on $\mathcal{B}_{o,d}$, respectively. Obviously, $E_{o,d} \subseteq E$ and $\mathcal{B}_{o,d} \subseteq \mathcal{B}$ hold for all $o,d \in V$. Then the scale of the MILP model can be reduced efficiently. Note that both set $E_{o,d}$ and $\mathcal{B}_{o,d}$ can be constructed easily in advance by generating all legal paths satisfying the detour ratio $\epsilon$ for each shipment. Meanwhile, we would not miss the optimal solution of the TBSP problem if every legal path is included.

\section{Sequential algorithm}
\label{sec: Sequential formulation}
\begin{figure}[h!]
  \caption{Framework of the sequential algorithm}\label{fig: Sequential algorithm}
  \centering
  \includegraphics[scale=0.82]{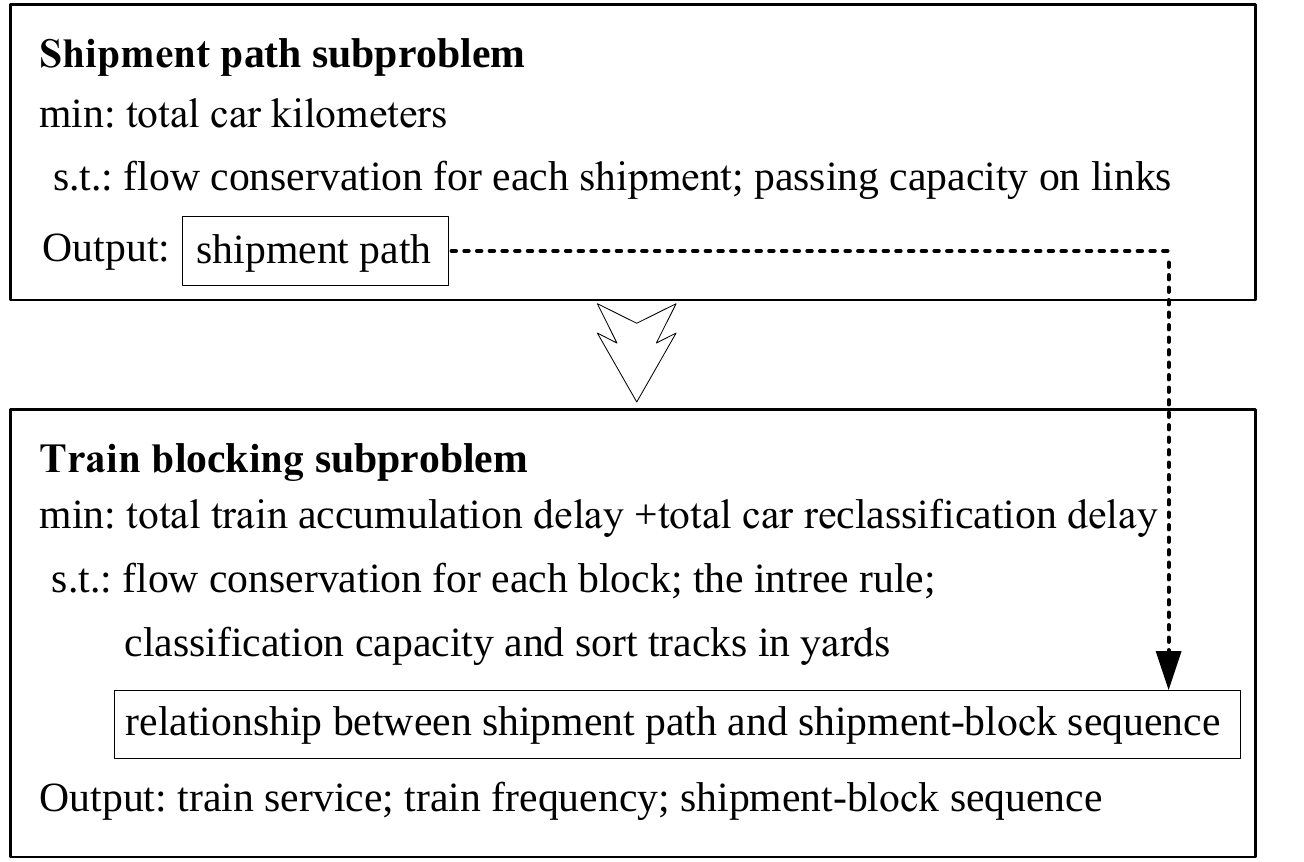}
\end{figure}

In this section, we decompose the TBSP problem into the shipment path subproblem and train blocking subproblem which are solved sequentially. The former seeks the path for each shipment to achieve the minimal total car kilometers. The latter identifies optimal train services and shipment-block sequences given the shipment paths, with the aim to minimize the total sum of train accumulation delay and car reclassification delay. The framework of the sequential algorithm is shown in Figure \ref{fig: Sequential algorithm}.

\subsection{Shipment path subproblem}
\label{sec: Routing subproblem}
The shipment path subproblem answers how to route all shipments over the physical network $G_1$ from their origin to their destination. Here, only the arc selection variables defined in Table \ref{tab: Decision Variables} are needed.
%($x^{o,d}_{i,j}, \forall o,d \in V, (i,j)\in E$)
The shipment path problem is modelled as a binary linear programming as follows.
\begin{equation}\small
\label{eq: Routing objective function}
\text{(Path)} \quad \min \sum_{o,d\in V}\sum_{(i,j)\in E}n_{o,d}l_{i,j}x_{i,j}^{o,d}
\end{equation}
\begin{align}\small
&\text{s.t.}&\sum_{o,d\in V }n_{o,d}x_{i,j}^{o,d} &\leq mf_{i,j}\alpha_{i,j} &&\forall (i,j) \in E \label{constraint: Routing link capacity}
\end{align}
\begin{equation*}\small
\eqref{constraint: car balance}, \eqref{constraint: detour}, \eqref{constraint: arc selection variables domain}.
\end{equation*}

The objective function \eqref{eq: Routing objective function} is to minimize the total car kilometers of shipping all demands over the physical network. Inequalities \eqref{constraint: Routing link capacity} ensure the usage of each link does not exceed its capacity which is measured by the total number of cars. Flow conservation equations \eqref{constraint: car balance} and detour inequations \eqref{constraint: detour} are also presented for each shipment. Formula \eqref{constraint: arc selection variables domain} specifies the domain of decision variables.
The Path model probably has $|V|^2|E|$ variables and $|V|^3+|V|^2+|E|$ constraints. We denote $R_{o,d}$ as the path of shipment from $o$ to $d$ computed by the shipment path subproblem.

\subsection{Train blocking subproblem}
\label{sec: Blocking problem}
The train blocking subproblem answers how to flow shipments on the service network $G_2$ along the paths from their origin to their destination given by the shipment path subproblem. It needs to makes the following decisions: which blocks %(pairs of yards)
are built, which shipments are grouped to which blocks. Train frequencies would be certainly determined by the total intensity of each block and the train size after solving the train blocking subproblem. Together with the consolidation selection and sort track usage decisions, there are four classes of variables in total whose notation are the same with those in Table \ref{tab: Decision Variables}. The train blocking subproblem is modelled as a mixed integer linear programming as follows.
\newpage
\begin{equation}\small
\label{eq: Blocking objective function}
\text{(Block)} \quad \min \sum_{p,q\in V}mc_{p,q}y_{p,q}
+ \sum_{o,d\in V}\sum_{\{p,q \in V | q\neq d\}}n_{o,d}t_qu^{o,d}_{p,q}
\end{equation}
\begin{align}\small
&\text{s.t.}&u_{p,q}^{o,d} &\leq y_{p,q}I(R_{p,q}\subseteq R_{o,d}) &&\forall o, d,p,q \in V \label{constraint: blocking car-blcok cross}
\end{align}
\begin{equation*}\small
  \eqref{constraint: block balance}-\eqref{constraint: exclusive rule 2}, \eqref{constraint: block design variables domain}, \eqref{constraint: car distribution variables domain}-\eqref{constraint: track usage variables domain}.
\end{equation*}

The objective function \eqref{eq: Blocking objective function} is to minimize the total sum of train accumulation delay and car reclassification delay. In constraint \eqref{constraint: blocking car-blcok cross}, the indicator function $I(R_{p,q}\subseteq R_{o,d})=1$ if $R_{p,q}\subseteq R_{o,d}$, 0 otherwise. There are two preconditions when shipment $o \rightarrow d$ is consolidated into block $p\rightarrow q$. One is that block $B_{p,q}$ has been provided, and another is the shorter path $R_{p,q}$ is completely contained in the longer path $R_{o,d}$. These two preconditions correspond to constraints \eqref{constraint: car-blcok cross} and \eqref{constraint: path consistant} in the integrated model, respectively. Besides, constraints \eqref{constraint: block balance}-\eqref{constraint: exclusive rule 2} are provided again which have been explained in subsection \ref{sec: MIQP formulation}. Formulas \eqref{constraint: block design variables domain} and \eqref{constraint: car distribution variables domain}-\eqref{constraint: track usage variables domain} specify the domain of decision variables. The Block model approximately has $|V|^4+|V|^3+2|V|^2$ variables and $2|V|^4+|V|^3+3|V|^2+2|V|$ constraints.

\subsection{Comparison between integrated model and sequential models}
\label{sec: Formulation comparison}
%In the shipment path subproblem, each shipment cannot be split and use multiple paths. Therefore, this subproblem could be viewed as an integer multi-commodity flow problem %\citep{Barnhart2000b}
%in which shipments play the act of commodities. On the other hand, the train blocking subproblem could be considered as the well-known service network design problem
%%\citep{Crainic2000}
%added with the intree constraint. Of course, both two subproblems are NP-hard.
In this subsection, we make a comparison between the integrated model and the two subproblem models in the sequential algorithm in terms of objectives, variables and constraints.
\begin{itemize}
  \item The objective function \eqref{eq: objective function} of the integrated model is the weighted sum of the one in the Path model and Block model, i.e. weighted sum of \eqref{eq: Routing objective function} and \eqref{eq: Blocking objective function}.
  \item Train frequency variables are only presented in constraints \eqref{constraint: link capacity} and \eqref{constraint: transport capacity} in the integrated model. Under the sequential framework, another version of link capacity \eqref{constraint: Routing link capacity} is provided in the Path model. Meanwhile, train frequencies are not independent decisions any more in the Block model.
  \item Constraints \eqref{constraint: car-blcok cross} and \eqref{constraint: path consistant} are comparable with \eqref{constraint: blocking car-blcok cross} to describe the logic relationship among block sequences, train services, and shipment paths.
  \item There are no longer quadric terms in either sequential models so that auxiliary variables \eqref{constraint: auxiliary variables domain 1} and constraints \eqref{constraint: auxiliary constraint 1}-\eqref{constraint: auxiliary constraint 3} in the integrated model are not required in the sequential algorithm.
\end{itemize}
Thus, compared with the integrated model, the sequential methodology has not only a far less model size but also a simpler model structure, which allows us to solve large-scale instances efficiently.

\section{Computational results}
\label{sec: Computational Results}
We now report the computational results of our proposed approaches to the three given data sets. The model and algorithm are coded in MATLAB 2016a. ILOG CPLEX 12.6 is used as the underlying optimization solver. And, we adopt the default parameter setting except the maximum runtime for the CPLEX. The MILP model and reduced MILP model are conducted on a GPU workstation with Inter Xeon E5 processor and 128 GB of RAM. The sequential algorithm is implemented on a Thinkpad T480 laptop with Intel Quad Core i7 processor and 8 GB of RAM.

\begin{figure}[h!]
  \caption{Number of feasible paths by shipments}\label{fig: Number of k shortest paths}
  \centering
  \includegraphics[scale=0.39]{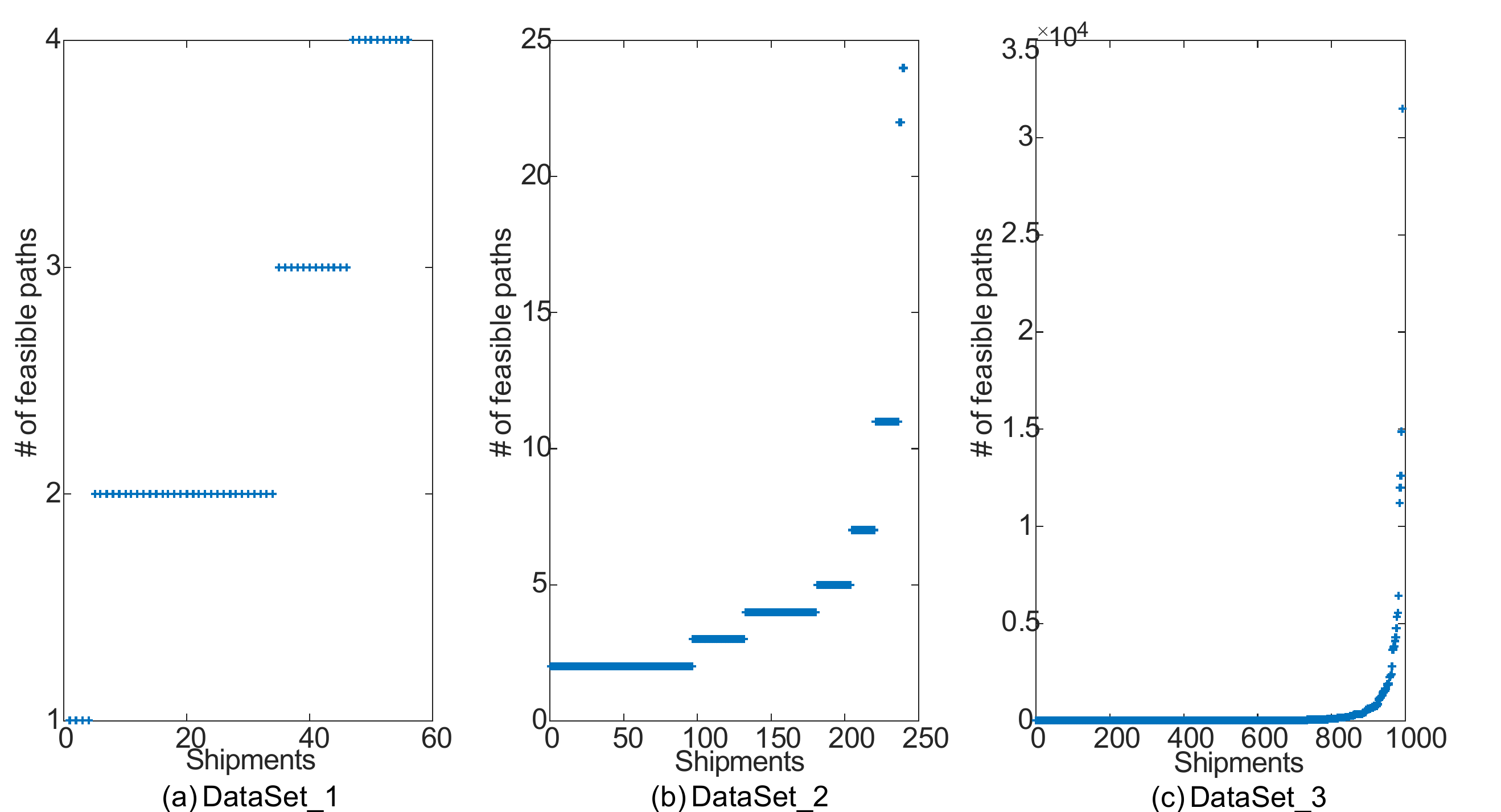}
\end{figure}

%\subsection{Comparison of model scale}
%We first compare the model scale of the MILP model, reduced MILP model and sequential algorithm on the tested data sets.
To implement the reduced model, we need to generate all feasible paths subject to the detour constraints \eqref{constraint: detour}. We apply the Yen's algorithm %\citep{Yen1971}
to produce all feasible paths under the given detour ratio for each shipment.
The number of feasible paths in the three data sets are displayed in Figure \ref{fig: Number of k shortest paths}.
%In each subfigure, the horizontal and vertical axis represent the index of shipments and the number of feasible paths of shipments, respectively.
%No demand is up to the preset value and the overwhelming majority are far from the threshold. As seen in this figure,the majority of shipments do not have too many feasible paths. Taking DataSet\_3 for example, about 93\% of shipments have less than 1000 feasible paths.
The scale of the MILP model, reduced MILP model and sequential algorithm are presented in Table \ref{tab: Comparison of model scale}. As observed in this table, the model scale of the reduced MILP model and sequential algorithm are far less than those of the MILP model. The model scale of the sequential algorithm would be reduced considerably when the network grows. Therefore, the storage requirements and computation effort are greatly reduced.
\begin{table}[ht]\small
\begin{center}
\caption{Comparison of model scale}
\label{tab: Comparison of model scale}
\begin{tabular}{llrrrrr}
\hline
\multirow{2}{*}{DataSet} & \multirow{2}{*}{Items} & \multirow{2}{*}{MILP} & \multirow{2}{*}{Reduced MILP}&\multicolumn{3}{c}{Sequential}\\
\cline{5-7}
& & & &Path & Block & Path+Block \\
\hline
\multirow{4}{*}{DataSet\_1}
&\# of variables&5,768 &1,712 &1,008 &452 &1,460\\
&\# of constraints&66,962 &5,640 &522 &800 &1,322\\
&\# of non-zero elements &\multirow{2}{*}{212,464} &\multirow{2}{*}{17,931} &\multirow{2}{*}{4,032} &\multirow{2}{*}{2,268} &\multirow{2}{*}{6,300}\\
&\textcolor{white}{\#} in coefficient matrix \\
\hline
\multirow{4}{*}{DataSet\_2}
&\# of variables&85,200 &11,307 &11,520 &2,480 &14,000\\
&\# of constraints&2,923,536 &42,888 &4,128 &4,352 &8,480\\
&\# of non-zero elements &\multirow{2}{*}{9,001,920} &\multirow{2}{*}{140,714} &\multirow{2}{*}{46,080} &\multirow{2}{*}{13,600} &\multirow{2}{*}{59,680}\\
&\textcolor{white}{\#} in coefficient matrix \\
\hline
\multirow{4}{*}{DataSet\_3}
&\# of variables&1,272,736 &151,800 &126,976 &15,739 &142,715\\
&\# of constraints&128,377,920 &1,465,843 &32,864 &27,968 &60,832\\
&\# of non-zero elements &\multirow{2}{*}{389,153,664} &\multirow{2}{*}{4,716,441} &\multirow{2}{*}{507,904} &\multirow{2}{*}{98,179} &\multirow{2}{*}{606,083 }\\
&\textcolor{white}{\#} in coefficient matrix \\
\hline
\end{tabular}
\end{center}
\end{table}

%\subsection{Computational results}
We set a maximum runtime of 12 hours on a GPU workstation for both the MILP model and reduced MILP model. Meanwhile, only 30 minutes on a laptop are allowed for the Path model and Block model of the sequential algorithm, respectively. The computational results of the three approaches on the data sets are listed in Table \ref{tab: Comparison of solution performance}. As shown in this table, both the MILP model and reduced MILP model quickly solve DataSet\_1 to optimality. The MILP model cannot find feasible solutions of DataSet\_2 and DataSet\_3 within the limited time. The reduced MILP model could provide strong lower bounds for the three data sets. The detour constraints on shipment paths is the starting point of the reduced reformulation. In Appendix \ref{app: Effect of detour constraints}, we investigate the solution quality of the TBSP problem by relaxation of the detour constraints. The sequential algorithm find a competitive solution with 0.04\% gap within 0.26 seconds for DataSet\_1, 0.42\% gap within 4.53 seconds for DataSet\_2 and 1.54\% gap within 0.58 hours for DataSet\_3, respectively. Thus, the sequential algorithm can solve the TBSP problem effectively within reasonable computation time.
\begin{table}[ht]\small
\centering
\begin{threeparttable}
\caption{Comparison of solution performance}
\label{tab: Comparison of solution performance}
\begin{tabular}{llrrr}
\hline
DataSet & Items & MILP& Reduced MILP &Sequential\\
\hline
\multirow{4}{*}{DataSet\_1}
&Lower bound (car-hour)&236,135 &236,135 &-- \\
&Upper bound (car-hour)&236,135 &236,135 &236,233 \\
&Gap$^{\rm a}$&0.00\%&0.00\%&0.04\%\\
&Run time  &9.44 s&8.48 s&0.26 s\\
\hline
\multirow{4}{*}{DataSet\_2}
&Lower bound (car-hour)&1,359,700 &1,362,565 &-- \\
&Upper bound (car-hour)&-- &1,364,772 &1,368,365 \\
&Gap&--&0.16\%&0.42\%\\
&Run time &12 h&6 h$^{\rm b}$&4.53 s\\
 \hline
\multirow{4}{*}{DataSet\_3}
&Lower bound (car-hour)&-- &5,090,775 &--\\
&Upper bound (car-hour)&-- &-- &5,170,224\\
&Gap&--&--&1.54\%\\
&Run time&3 h$^{\rm b}$&12 h&0.58 h\\
\hline
\end{tabular}
\begin{tablenotes}
\item[a] Gap=(Upper bound - Lower bound)/Upper bound *100\%.
\item[b] Out of memory.
\end{tablenotes}
\end{threeparttable}
\end{table}

\section{Conclusions}
\label{sec: Conclusions}
We address the integrated train blocking and shipment path optimization for tonnage-based operating railways. We firstly develop a mixed integer linear programming model that incorporates the intree rule in an integrated manner. Some model reduction technologies are discussed to cut down the model scale effectively. We further develop a sequential algorithm that decomposes the original problem into the shipment path subproblem and the train blocking subproblem which are solved sequentially. In the first stage, the shipment path subproblem finds out the least transportation cost for the shipment paths while coordinating all lines' passing capacity. In the second stage, the train blocking subproblem determines the train services, train frequencies and shipment-block sequences, with respect to the limited classification capacity and sort tracks in each yard as well as the intree rule. Computational results of the three data sets show that, the reduced MILP could provide strong lower bounds for the large-scale instances, and the sequential algorithm can derive competitive feasible solutions ($\thicksim 1\%$ gap) within reasonable computation time ($<1$ h).

%\section*{Acknowledgment}
%The research was supported by the National Natural Science Foundation of China (No. 51505309, 61603318), the Fundamental Research Funds for the Central Universities (No. 2682016CX118).

\bibliographystyle{plainnat}
\bibliography{biblio}\small

\begin{appendices}
\section{Effect of detour constraints}
\label{app: Effect of detour constraints}
In this appendix, we compare the solution quality of the TBSP problem under the two scenarios with and without the detour constraints \eqref{constraint: detour}. The three data sets under such two scenarios are both solved by the sequential algorithm on a laptop, and 30 minutes are allowed for the Path model and Block model. Table \ref{tab: Detour constraints relaxation} summarizes the comparison results. As indicated in this table, the detour constraints do not impact the solution quality of DataSet\_1.
For DataSet\_2 and DataSet\_3, compared with the scenario with \eqref{constraint: detour}, the relaxation one has a little increase of the car mile and car reclassification delay but a larger decrease in the train accumulation delay. This results in lower total costs for DataSet\_2 and DataSet\_3.
%For DataSet\_2 and DataSet\_3, compared with the scenario with \eqref{constraint: detour}, the scenario without \eqref{constraint: detour} has some increase in the car mile and train accumulation delay and some increase in the car reclassification delay. The relaxation of the detour constraints results in different consequences on the total cost of DataSet\_2 and DataSet\_3. For DataSet\_2, the relaxation model could find better solutions because of the larger and under-control feasible regions. However, for DataSet\_3, the relaxation model could hardly find better solutions because of the larger but out-of-control feasible regions under the limitation of run time.
Overall, we conclude that the specified detour ratio has no notable impact on the solution quality of the three given data sets.
\begin{table}[ht]\small
\begin{center}
\caption{Comparison of detour constraints}
\label{tab: Detour constraints relaxation}
\begin{tabular}{llrrr}
\hline
DataSet & Items & With \eqref{constraint: detour} & Without \eqref{constraint: detour} & Relative deviation\\
 \hline
\multirow{5}{*}{DataSet\_1}
&Car mile (car-km)&2,119,046 &2,119,046 &0.00\%\\
&Accumulation (car-hour)&19,700 &19,700 &0.00\%\\
&Classification cost (car-hour)&4,629 &4,629 &0.00\%\\
&Total cost (car-hour)&236,233 &236,233 &0.00\%\\
&Run time  (second)&0.26 &0.39 &50.00\%\\
  \hline
\multirow{5}{*}{DataSet\_2}
&Car mile (car-km)&12,537,081 &12,537,163 &0.00\%\\
&Accumulation (car-hour)&99,072 &98,632 &-0.44\%\\
&Classification cost (car-hour)&15,585 &15,953 &2.36\%\\
&Total cost (car-hour)&1,368,365 &1,368,301 &0.00\%\\
&Run time  (second)&4.53 &7.02 &54.97\%\\
  \hline
\multirow{5}{*}{DataSet\_3}
&Car mile (car-km)&47,999,769 &47,999,801 &0.00\%\\
&Accumulation (car-hour)&249,125 &245,330 &-1.52\%\\
&Classification cost (car-hour)&121,123 &123,457 &1.93\%\\
&Total cost (car-hour)&5,170,224 &5,168,767 &-0.03\%\\
&Run time  (second)&2088 &2285 &9.43\%\\
\hline
\end{tabular}
\end{center}
\end{table}

\end{appendices}

\end{document}